\documentclass[a4paper, 11pt] {article}
\usepackage[cp1251]{inputenc}
\usepackage[russian, english]{babel}
\usepackage{amsfonts}
\usepackage{mathtext}
\usepackage{cite}

\usepackage{graphicx}

\usepackage{verbatim}
\usepackage{amsmath}
\usepackage{amsthm}
\usepackage{amssymb}
\usepackage{delarray}
\usepackage{mathrsfs}
\usepackage{xcolor}

\begin{document}

\thispagestyle{empty}

\begin{center}
{\bf\Large Uniform stability of recovering the Sturm-Liouville operators with frozen argument}
\end{center}

\begin{center}
{\bf\large Maria Kuznetsova\footnote{Department of Mathematics, Saratov State University, e-mail: {\it kuznetsovama@info.sgu.ru}}}
\end{center}

\noindent{\bf Abstract.} 
In the paper, we study the problem of recovering the Sturm--Liouville operator with frozen argument from its spectrum and additional data. 
For this inverse problem, we establish a substantial property of the uniform stability, which consists in that the potential depends Lipschitz continuously on the input data. 

\smallskip
\noindent {\it Key words}: inverse spectral problems; frozen argument; Sturm--Liouville operators; uniform stability

\smallskip
\noindent {\it 2010 Mathematics Subject Classification}: 34K29, 34A55
\\

{\large\bf 1. Introduction}
\\

Inverse spectral problems consist in recovering operators from their spectral characteristics. The most complete results in this direction concern the classical  differential operators  (see~\cite{B,yurko,march,levitan}).
Here, we consider the Sturm--Liouville operator with frozen argument, which, unlikely the classical ones, is nonlocal. 
By this reason, investigation of inverse spectral problems for the operator with frozen argument requires nontraditional methods.

We study the recovery of a complex-valued potential $q \in L_2(0, \pi)$ from the spectrum $\{ \lambda_n \}_{n \in \mathbb N}$ of the boundary value problem ${\cal B}_{a,\alpha,\beta}(q):$
\begin{equation}
\ell y := -y''(x) + q(x) y(a) = \lambda y(x), \quad x \in (0, \pi), 
\label{SL-equation}
\end{equation}
\begin{equation} \label{bc}
y^{(\alpha)}(0) = y^{(\beta)}(\pi) = 0,
\end{equation}
where $a \in [0, \pi]$ and $\alpha, \beta \in \{ 0, 1 \}.$  Expression $\ell y$ determines the Sturm--Liouville operator with frozen argument.  At the same time, equation~\eqref{SL-equation} belongs to the class of the so-called
loaded differential equations \cite{nahushev, kraal, lomov, pol}, which have application in mathematics, physics, mathematical
biology, etc. 

The aim of this paper is to obtain the stability of the inverse problem by the spectrum and some additional data. 
This property guaranties that small perturbations
of any fixed input data lead to small
perturbations of the solution and has a usage in the development of numerical algorithms, see~\cite{delay num, ignat, bon-gaidel}.
Here, we prove the stronger uniform stability involving uniform estimates.
Earlier, various authors addressed the uniform stability (see~\cite{sav,hryn,hryn2,B1,B2,but-djur,graph-delay}), that allowed them to better comprehend the nature of studied inverse problems.  

Various aspects of inverse spectral problems for operators with frozen argument were studied in~\cite{albeverio, nizh09, nizh10, BBV, BV,
BK, wang, but-hu, tsai, my-froz, bon-disc}. 
In particular, the inverse problem by the spectrum $\{ \lambda_n \}_{n \in \mathbb N}$ was investigated in~\cite{BBV, BV, BK, wang,
tsai}, wherein the cases of rational and irrational $a/\pi$ were treated apart by different techniques. 
In the recent work~\cite{AML}, there was offered a new unified approach to the both cases, that let us to prove the spectrum characterization:
An arbitrary sequence $\{ \lambda_n \}_{n \in \mathbb N}$ is the spectrum  of some boundary value problem ${\cal B}_{a,\alpha,\beta}(q)$ if and only if
\begin{equation} \label{main asymp}
\lambda_n = \rho_n^2, \quad \rho_n = z_n + \frac{\nu_n}{n} \varphi_\alpha(a z_n), \quad n \in {\mathbb N}, \quad \{\nu_n \}_{n \in {\mathbb N}} \in \ell_2, 
\end{equation}
where we denoted
$$z_n := n - \frac{\alpha+\beta}{2}, \quad \varphi_\alpha(x) := \left\{ \begin{array}{cc}
\sin x, & \alpha = 0, \\
\cos x, & \alpha = 1.
\end{array}\right. $$
We also mention another general approach in~\cite{DobHry} that allows to obtain the spectrum characterization. 

Formula~\eqref{main asymp} gives peculiar asymptotics with the multiplier $\varphi_\alpha(a z_n)$ in the $n$-th residual term, $n \in \mathbb N.$
If $a/\pi$ is irrational, the multiplier never turns $0$ though can be arbitrarily small by modulus. For rational $a/\pi,$ the values of $\varphi_\alpha(a z_n)$ are cyclic and may include $0.$
Put $\Omega=\{n\in \mathbb N \colon  \varphi_\alpha(az_n) =0\}$ and $\overline{\Omega} = \mathbb{N}\setminus\Omega.$ 
Formula~\eqref{main asymp} yields that the spectrum part $\{ \lambda_n \}_{n \in \Omega} $ is the same for any $q,$ i.e. this part is uninformative.
In~\cite{AML}, it was proved that for the unique recovering of $q,$ one should additionally set the part $\{ \xi_n\}_{n \in \Omega} \in \ell_2$ of the coefficients $\xi_n = \int_0^\pi \varphi_\alpha(z_n t) q(t) \, dt,$ $n \in \mathbb N.$  Thus, we have to consider the following inverse problem.
\medskip

\noindent{\bf Inverse problem~1.}
Given $\{ \lambda_n \}_{n \in \overline{\Omega}}$ and $\{ \xi_n\}_{n \in \Omega},$ recover $q \in L_2(0, \pi).$
\medskip

Let $\| \cdot \|$ be the classical norm in the space $\ell_2$  and
$$
\big\| \{ x_n\}_{n \in \overline{\Omega}}\big\|_{\bf a} := 
\left\| \big\{ x_n \varphi^{-1}_\alpha(az_n)\big\}_{n \in \overline{\Omega}}\right\|.
$$
Then, the value $\big\| \{ n(\rho_n - z_n)\}_{n \in \overline{\Omega}}\big\|_{\bf a}$ is finite provided~\eqref{main asymp} holds.

We prove the following theorem on the uniform stability of Inverse problem~1. 
\medskip

\noindent{\bf Theorem~1.} {\it
Fix $r > 0$ and consider two boundary value problems ${\cal B}_{a,\alpha,\beta}(q)$ and ${\cal B}_{a,\alpha,\beta}(\tilde q)$
having spectra $\{ \rho^2_n \}_{n \in \mathbb N}$ and $\{\tilde\rho^2_n \}_{n\in \mathbb N},$ respectively. If the condition 
 \begin{equation} \label{bounded}
\big\|\{ n(\rho_n-z_n) \}_{n\in \overline{\Omega}}\big\|_{\bf a} \le r, \quad \big\|\{  n(\tilde\rho_n-z_n) \}_{n\in \overline{\Omega}}\big\|_{\bf a} \le r
\end{equation}
is fulfilled, then we have the inequality
\begin{equation} \label{result}
\| q - \tilde q\|_{L_2(0, \pi)} \le C_r \Xi + \sqrt{\frac{2}{\pi}}\big\| \{\xi_n - \tilde \xi_n\}_{n \in \Omega}\big\|, \quad \Xi := \big\| \{n(\rho_n - \tilde\rho_n)\}_{n \in \overline{\Omega}}\big\|_{\bf a},
\end{equation}
where $\xi_n = \int_0^\pi q(t) \varphi_\alpha(z_n t) \, dt,$ $\tilde\xi_n = \int_0^\pi \tilde q(t) \varphi_\alpha(z_n t) \, dt,$ 
and the constant $C_r$ depends only on $r.$ 
}
\medskip

The statement of Theorem~1 means that the potential $q$ depends Lipschitz continuously on the input data of Inverse problem~1 from the sets determined by the conditions~\eqref{main asymp}, ${\big\| \{ n(\rho_n - z_n)\}_{n \in \overline{\Omega}}\big\|_{\bf a} \le r,}$ and $\{ \xi_n\}_{n \in \Omega} \in \ell_2.$

 The uniform stability of recovering non-local operators of other types was obtained in~\cite{but-djur,B1,B2,graph-delay}. 
The approach to nonlocal operators developed there emloyed the uniform stability of reconstructing the characteristic functions from the spectra. We follow the same strategy. Compared to the works~\cite{but-djur,B1,B2,graph-delay}, here the specifics is in the value $\Xi$ used for measuring the distance between spectra. Its defintion involves division on $\varphi_\alpha(z_na),$  being potentially arbitrarily close to $0.$ This leads to the need in special estimates for the differences of the characteristic functions in the proof of Theorem~1.
\\

{\large\bf 2. Preliminaries}
\\

Consider the function $\Delta(\lambda)$ entire in $\lambda:$ 
\begin{multline*} 
\Delta(\lambda) = (-1)^{\alpha(1-\beta)}\frac{\varphi_{|\beta-\alpha|}(\pi \rho)}{\rho^{1-\alpha-\beta}} + \\
+ (-1)^\alpha\Bigg(\frac{\varphi_{\beta}((\pi-a) \rho)}{\rho^{2-\alpha-\beta}} \int_0^a q(t) \varphi_{\alpha}(t\rho) \, dt
+\frac{\varphi_{\alpha}(a\rho)}{\rho^{2-\alpha-\beta}} \int_a^\pi q(t) \varphi_{\beta}((\pi - t)\rho) \, dt\Bigg),
\end{multline*}
where and below $\lambda=\rho^2.$ In~\cite{AML,wang}, it was proved that $\Delta(\lambda)$ is the characteristic function of the boundary value problem ${\cal B}_{a, \alpha,\beta}(q),$ i.e. a number $\lambda_n$ is an eigenvalue if and and only if $\Delta(\lambda_n)=0.$ The sequence of the eigenvalues $\{ \lambda_n \}_{n \in \mathbb N}$ of the boundary value problem ${\cal B}_{a, \alpha,\beta}(q)$ taken with the account of algebraic multiplicities is called the spectrum.

For the characteristic function, we obtain the representation
\begin{equation}\label{Delta}
\Delta(\lambda) = (-1)^{\alpha(1-\beta)}\frac{\varphi_{|\beta-\alpha|}(\pi \rho)}{\rho^{1-\alpha-\beta}} + \frac{1}{\rho^{2-\alpha-\beta}}\int_0^\pi \varphi_{1-|\beta-\alpha|}(t \rho) \, W(t)\,dt, \quad W \in L_2(0, \pi),
\end{equation}
wherein $\int_0^\pi W(t) \, dt = 0$ if $\alpha=\beta=0.$
The proof of formula~\eqref{Delta} was given in~\cite{BV,BBV} for rational  $a/\pi;$ for irrational $a/\pi,$ the computations are analogous.
It is known that $\Delta(\lambda)$ and, in turn, $W(t)$ are uniquely recovered from the spectrum $\{ \lambda_n \}_{n \in \mathbb N}$ (see~\cite{BK, BV, BBV,wang}). Further, we have to obtain the uniform stability of recovering $W(t)$ 
with respect to the values ${\rho_n = \sqrt{\lambda_n}},$ $n \in \mathbb N,$ where $\arg \rho_n \in (-\frac\pi2, \frac\pi2].$

 Consider two spectra $\{\lambda_n\}_{n \in \mathbb N}$ and $\{ \tilde \lambda_n\}_{n \in \mathbb N}.$
We agree that if a certain symbol $\gamma$ denotes an object related to  $\{ \lambda_n \}_{n \in \mathbb N},$ then this symbol with tilde
$\tilde\gamma$ will denote the analogous object related to $\{ \tilde \lambda_n\}_{n \in \mathbb N}.$ For briefness, we designate $\hat\gamma := \gamma - \tilde \gamma.$ We also use one and the same notation $C_r$ for various constants depending only on $r.$
\medskip

\noindent{\bf Lemma~1.} {\it
Let 
\begin{equation} \label{l2 bounded}
\big\| \{ n (\rho_n - z_n) \}_{n \in \mathbb N}\big\| \le r, \quad \big\| \{ n (\tilde \rho_n - z_n)\}_{n \in \mathbb N}\big\| \le r,
\end{equation} 
 where $r>0$ is fixed. Then, the following inequalities hold:}
\begin{equation} \| W\|_{L_2(0, \pi)} \le C_r, \quad \| \tilde W\|_{L_2(0, \pi)} \le C_r; \quad 
\| \hat W\|_{L_2(0, \pi)} \le C_r \big\| \{ n \hat \rho_n\}_{n \in \mathbb N} \big\|.
\label{W}
\end{equation}

Lemma~1 is obtained as a consequence of  ~\cite[Theorem~7]{MN}, since $\Delta(\rho^2)$ is a particular case of sine-type functions studied in~\cite{MN}.
For the proof of Theorem~1, we also need the following lemma.
\medskip

\noindent{\bf Lemma~2.} {\it 
Let $M>0.$  If $\| W\|_{L_2(0, \pi)} \le M$ and $\Delta(\lambda)$ is given by~\eqref{Delta}, then}
\begin{equation} \label{est}
| \Delta(\rho^2)| \le C_M n^{\alpha+\beta-1} \text{ for } |\rho - z_n| \le M, \; n \in \mathbb N.
\end{equation}
\begin{proof}
Consider the case $\alpha=\beta=0$ (the other cases are proceeded analogously). From~\eqref{Delta} we have 
\begin{equation} \label{Delta00}
\Delta(\rho^2) = \frac{\sin \rho \pi}{\rho} + \int_0^\pi \frac{\cos \rho t}{\rho^2} W(t) \, dt, \quad \int_0^\pi W(t) \,dt = 0,\end{equation}
and $z_n = n.$ 
Let $|\rho - n| \le M$ for some $n \in \mathbb N.$  

If $n - M > 2,$ then for $|\rho - n| \le M,$ we have
$$|\rho|^{-1} \le (n - M)^{-1} \le C_M n^{-1}; \quad |\sin \rho \pi| < C_M; \quad |\cos \rho t| < C_M, \; t \in [0, \pi].$$ 
Substituting these estimates into~\eqref{Delta00}, we obtain
$|\Delta(\rho^2)| \le C_M n^{-1} + C_M n^{-2} \| W\|_{L_2(0, \pi)} \le C_M n^{-1},$
and~\eqref{est} is proved for $n - M > 2.$

If $n - M \le 2,$ then for $|\rho - n| \le M,$ the inequality $|\rho| \le 2(M+1)$ holds.
Consider the circle $R_M = \big\{ \lambda \in \mathbb C \colon |\lambda| \le 4(M+1)^2\big\}$ in the $\lambda$-plane. For $\lambda=\rho^2 \in R_M,$ we have the estimates
$$|\rho \sin \rho \pi| \le C_M, \quad |\cos \rho t|\le C_M, \; t \in [0, \pi],$$ 
and, consequently,
$|\lambda \Delta(\lambda)| \le C_M.$
Applying the Schwartz lemma to the function $\lambda \Delta(\lambda)$ in $R_M,$ we arrive at $|\Delta(\lambda)| \le C_M$ for $\lambda \in R_M.$ Since $n \le M+2,$ this yields
$|\Delta(\rho^2)| \le C_M n^{-1}$ for $|\rho| \le 2(M+1),$ and for $|\rho - n| \le M.$ The lemma is proved.
\end{proof}

{\large\bf 3. The proof of Theorem~1}
\\

Note that the system  $\{\varphi_\alpha(z_n t)\}_{n \in \mathbb N}$ becomes an orthogonal basis after dividing by $\sqrt{\pi/2}$ (if $\alpha = \beta=1,$ the first function should be divided by $\sqrt\pi$).
Then, we have the inequality
\begin{equation} \label{hat q}
\| \hat q \|_{L_2(0, \pi)} \le \sqrt\frac{2}{\pi} \big\|\{ \hat \xi_n\}_{n \in \mathbb N}\big\| \le 
\sqrt\frac{2}{\pi} \Big(\big\|\{ \hat \xi_n\}_{n \in \Omega}\big\|+\big\|\{ \hat \xi_n\}_{n \in \overline{\Omega}}\big\|\Big).
\end{equation}
 
In~\cite{AML}, it was proved that $\xi_n = (-1)^{n+1-\alpha} z_n^{2-\alpha-\beta} \Delta(z_n^2)\varphi^{-1}_\alpha(z_n a),$ $n \in \overline{\Omega},$
which yields 
\begin{equation*} 
\hat \xi_n =  (-1)^{n+1-\alpha} z_n^{2-\alpha-\beta} \frac{\hat \Delta(z_n^2)}{\varphi_\alpha(z_n a)}, \quad n \in \overline{\Omega}.
\end{equation*}
Thus, we have to estimate 
\begin{equation}\label{main relation}
\big\|\{ \hat \xi_n\}_{n \in \overline\Omega}\big\|=\big\| \{ \hat \Delta(z_n^2)z_n^{2-\alpha-\beta}\}_{n \in \overline\Omega}\big\|_{\bf a}.\end{equation}

 Further, for $n \in \overline\Omega,$ we consider each value $\hat \Delta(z_n^2)$ individually, representing it in the form
\begin{equation} \label{hat Delta}
\hat \Delta(z_n^2) = -\tilde d_n \hat \rho_n+ \hat d_n (z_n - \rho_n), \quad 
\tilde d_n := \lim_{\rho \to z_n} \tilde f(\rho), \quad \tilde f(\rho) := \frac{\tilde\Delta(\rho^2)}{\rho - \tilde\rho_n},
\end{equation}
while for $d_n$ and $f(\rho),$ we use the analogous notations without tilde.
Let us estimate $\tilde d_n$ and $\hat d_n$ apart with the Schwartz lemma.

a) It is clear that $\| \{ x_n \}_{n \in \overline{\Omega}} \| \le \| \{ x_n\}_{n \in \overline{\Omega}} \|_{\bf a}.$ Moreover, by~\eqref{main asymp}, we have
$$\big\| \{ n (\rho_n - z_n) \}_{n \in \mathbb N}\big\| = \big\| \{ n (\rho_n - z_n) \}_{n \in \overline{\Omega}}\big\| \le \big\| \{ n (\rho_n - z_n) \}_{n \in \overline{\Omega}}\big\|_{\bf a}.$$ 
Then, condition~\eqref{bounded} yields~\eqref{l2 bounded},
and by virtue of Lemma~1, inequalities~\eqref{W} hold.
Consider three following circles:
$$\tilde S_n = \{ z \in \mathbb C \colon |z - \tilde \rho_n| < 2r\}, \quad   S_n = \{ z \in \mathbb C \colon |z -  \rho_n| < 4r\}, \quad Z_n = \{ z \in \mathbb C \colon |z - z_n| < 6r\}.$$
Inequalities~\eqref{l2 bounded} yield $z_n \in \tilde S_n$ and $\tilde S_n \subset S_n \subset Z_n.$ 
By~\eqref{W} and Lemma~2, for the function $\tilde\Delta(\rho^2)$ in $Z_n,$ we have estimate~\eqref{est} with $M=\max(C_r, 6r).$ Applying the Schwartz lemma to the function $\tilde\Delta(\rho^2)$ in $\tilde S_n,$ we get
\begin{equation} \label{d est}
\big|\tilde f(\rho)\big| \le C_r n^{\alpha+\beta-1}, \; \rho \in \tilde S_n; \quad \big|\tilde d_n\big| = \Big|\lim_{\rho \to z_n}\tilde f(\rho)\Big| \le C_r n^{\alpha+\beta-1}.
\end{equation}
Analogously, we have 
\begin{equation} \label{f}
|f(\rho)| \le C_r n^{\alpha+\beta-1}, \quad \rho \in  S_n.
\end{equation}

b) Note that
$$\hat d_n = \lim_{\rho \to z_n} \frac{g(\rho)}{\rho - \tilde \rho_n}, \quad g(\rho) := \hat \Delta(\rho^2) + \hat \rho_n f(\rho), \quad g(\tilde \rho_n) = 0.$$

From formula~\eqref{Delta} it follows that
$\hat \Delta(\rho^2) = \rho^{\alpha+\beta-2} \int_0^\pi \varphi_\alpha(\rho t) \hat W(t) \, dt.$
Proceeding analogously to the proof of Lemma~2, using~\eqref{W}, we arrive at the inequality 
\begin{equation} \label{g1}
\big|\hat\Delta(\rho^2)\big| \le \frac{C_r}{n^{2 - \alpha-\beta}}\| \hat W\|_{L_2(0, \pi)} \le \frac{C_r}{n^{2 - \alpha-\beta}} \Xi, \quad \rho \in Z_n.\end{equation}

For estimating the second summand in $g(\rho),$ we use~\eqref{f} and that $|n \hat \rho_n|   \le \Xi:$
\begin{equation} \label{g2}
|\hat\rho_n f(\rho)| \le \frac{C_r}{n^{2-\alpha-\beta}} \Xi, \quad \rho \in S_n. 
\end{equation}
Combining~\eqref{g1} and~\eqref{g2}, we obtain $|g(\rho)| \le C_r n^{\alpha+\beta-2} \Xi$ in $\tilde S_n.$ By the Schwartz lemma,  
\begin{equation} \label{g}
\left|\frac{g(\rho)}{\rho - \tilde \rho_n}\right| \le C_r n^{\alpha+\beta-2} \Xi, \; \rho \in \tilde S_n; \quad \big|\hat d_n\big| \le  C_r n^{\alpha+\beta-2} \Xi.\end{equation}

Finally, unifying~\eqref{d est} and~\eqref{g} with~\eqref{hat Delta}, we have
$$\big|\hat \Delta(z_n^2)\big| \le \frac{C_r |\hat \rho_n|}{n^{1-\alpha-\beta}} + \frac{C_r \Xi |z_n - \rho_n|}{n^{2-\alpha-\beta}}, \quad n \in \overline\Omega.$$
Multiplying by $n^{2-\alpha-\beta}\varphi^{-1}_\alpha(z_n a)$ and taking $\ell_2$-norm, we arrive at
$$\big\| \{ \hat \Delta(z_n^2) z_n^{2-\alpha-\beta}\}_{n \in \overline\Omega}\big\|_{\bf a} \le C_r\Big(\big\| \{n\hat \rho_n\}_{n\in \overline\Omega}\big\|_{\bf a} + \Xi \big\| \{n (\rho_n-z_n)\}_{n \in \overline\Omega}\big\|_{\bf a}\Big) \le C_r \Xi.$$
Application of estimates~\eqref{hat q} and~\eqref{main relation} finishes the proof. $\qed$
\medskip

\noindent{\it Remark.} Analogously, one can study the stability of recovering $q$ involving estimates with $\{ \lambda_n\}_{n \in \overline{\Omega}}$ instead of $\{ \rho_n\}_{n \in \overline{\Omega}}.$ In particular, in the case $\alpha=\beta=0$ we proved the following theorem.
\medskip

\noindent{\bf Theorem 2.} {\it
Fix $r > 0$ and consider two boundary value problems ${\cal B}_{a,0,0}(q)$ and ${\cal B}_{a,0,0}(\tilde q)$ having the spectra $\{ \lambda_n \}_{n \in \mathbb N}$ and $\{ \tilde\lambda_n \}_{n \in \mathbb N},$ respectively.  If the condition 
$$ \big\|\{ \lambda_n - n^2 \}_{n \in \overline{\Omega}}\big\|_{\bf a} \le r, \quad \big\|\{ \tilde \lambda_n - n^2 \}_{n \in \overline{\Omega}}\big\|_{\bf a} \le r$$ 
is fulfilled, then we have the inequality
\begin{equation} \label{lambda stab}
\| \hat q\|_{L_2(0,\pi)} \le C_r \big\| \{\hat \lambda_n\}_{n \in \overline{\Omega}} \big\|_{\bf a} + \sqrt{\frac{2}{\pi}}\big\| \{\hat \xi_n\}_{n \in \Omega}\big\|.
\end{equation}
 }

Actually, this theorem does not follow from Theorem~1 without additional assumptions. 
For obtaining it, we adapted the proof of Theorem~1.
The most significant alteration is that we had to prove the analogue of Lemma~1 with  $\big\| \{ \hat \lambda_n\}_{n \in \mathbb N} \big\|$ instead of  $\big\| \{ n \hat \rho_n\}_{n \in \mathbb N} \big\|$ on the right side of~\eqref{W}.
For the case $\alpha=\beta=0,$ the required statement was proved in~\cite{B1}.

\smallskip

\noindent{\bf Acknowledgement.} This research was supported by a grant of the Russian Science Foundation No.~22-21-00509,
https://rscf.ru/project/22-21-00509/. The author thanks Sergey Buterin for valuable comments.


\begin{thebibliography}{99}

\bibitem{B}
Borg G. Eine Umkehrung der Sturm--Liouvilleschen Eigenwertaufgabe, Acta Math. 78 (1946), 1--96.

\bibitem{march}
Marchenko V.A. Sturm--Liouville Operators and Their Applications, Naukova Dumka, Kiev, 1977; English transl., Birkh\"{a}user, 1986.

\bibitem{levitan}
Levitan B.M. Inverse Sturm--Liouville Problems, Nauka, Moscow, 1984; English transl., VNU Sci.Press, Utrecht, 1987.

\bibitem{yurko}
Freiling G., Yurko V.A. Inverse Sturm--Liouville Problems and Their Applications, NOVA Science Publishers, New York, 2001.

\bibitem{kraal}
Krall  A.M. The development of general differential and general differential-boundary systems, Rocky Mountain J. Math. 5 (1975), 493--542.

\bibitem{nahushev}
Nakhushev A.M. Loaded Equations and Their Applications, Nauka, Moscow, 2012.

\bibitem{lomov}
Lomov I.S. Loaded differential operators: Convergence of spectral expansions, Differential Equations 50 (2014), no. 8, 1070--1079.


\bibitem{pol}
Polyakov D.M. Nonlocal perturbation of a periodic problem for a second-order differential operator, Differential Equations 57 (2021),
11--18.

\bibitem{ignat} Ignatiev M., Yurko V. Numerical methods for solving inverse Sturm--Liouville problems, Result. Math. 52 (2008), 63--74 . 

\bibitem{delay num} Bondarenko N., Buterin S.
Numerical solution and stability of the inverse spectral problem for a convolution integro-differential operator,
 Commun. Nonlinear Sci. Numer. Simul. 89 (2020), 105298.

\bibitem{bon-gaidel} Bondarenko N.P., Gaidel A.V.  Solvability and stability of the inverse problem for the quadratic
differential pencil, Mathematics 9 (2021), no. 20, article 2617.

\bibitem{sav} Savchuk A.M., Shkalikov A.A. Inverse problems for Sturm--Liouville operators with potentials in Sobolev spaces: Uniform stability, Funct Anal Its Appl 44 (2010), 270--285. 

\bibitem{hryn} Hryniv R.O. Analyticity and uniform stability in the inverse spectral problem for Dirac operators, J. Math. Phys. 52 (2011), article~063513.

\bibitem{hryn2} Hryniv  R.O. Analyticity and uniform stability in the inverse singular Sturm--Liouville spectral problem, Inverse Problems 27 (2011), article~065011.

\bibitem{B1} Buterin S. Uniform stability of the inverse spectral problem for a
convolution integro-differential operator, Appl. Math. Comput. 390 (2021),
article~125592.

\bibitem{B2}
Buterin S.A. Uniform full stability of recovering convolutional perturbation of the Sturm--Liouville operator from the spectrum, Journal of
Differential Equations 282 (2021), 67--103.

\bibitem{but-djur} Buterin S., Djuri\'c N., Inverse problems for Dirac operators with constant delay: uniqueness, characterization, uniform stability, Lobachevskii J. Math. 43 (2022), no.~6, 1492--1501.

\bibitem{graph-delay} Buterin S. Functional-differential operators on geometrical graphs with global delay and inverse spectral problems, Results Math. 78 (2023), article~79.


\bibitem{albeverio}
Albeverio S., Hryniv R.O., Nizhnik L.P. Inverse spectral problems for non-local Sturm--Liouville operators, Inverse Problems 23 (2007), no.
2, 523--535.

\bibitem{nizh09}
Nizhnik L.P. Inverse eigenvalue problems for nonlocal Sturm--Liouville operators, Meth. Func. Anal. Top. 15 (2009), no.1, 41--47.

\bibitem{nizh10}
Nizhnik L.P., Inverse nonlocal Sturm-Liouville problem, Inverse Problems 26 (2010), no. 12, 125006.

\bibitem{BBV}
Bondarenko N.P., Buterin S.A., Vasiliev S.V. An inverse spectral problem for Sturm--Liouville operators with frozen argument, J. Math.
Anal. Appl. (2019), no. 1, 1028--1041.

\bibitem{BV}
Buterin S.A., Vasiliev S.V. On recovering a Sturm-Liouville-type operator with the frozen argument rationally proportioned to the interval
length, J. Inv. Ill-posed Probl., 27 (2019), no. 3, 429--438.

\bibitem{BK}
Buterin S., Kuznetsova M. On the inverse problem for Sturm--Liouville-type operators with frozen argument: rational case, Comp. Appl. Math. 39 (2020), article 5. 


\bibitem{but-hu}
Buterin S., Hu Y.T. Inverse spectral problems for Hill-type operators with frozen argument, Anal. Math. Phys. 11 (2021), article 75. 

\bibitem{wang}
Wang Y.P, Zhang M., Zhao W., Wei X. Reconstruction for Sturm--Liouville operators with frozen argument for irrational cases, Applied
Mathematics Letters 111 (2021), 106590.

\bibitem{tsai}
 Tsai T.M., Liu H.F., Buterin S., Chen L.H, Shieh C.T. Sturm--Liouville-type operators with frozen argument and Chebyshev polynomials, Math. Meth. Appl. Sci. 45 (2022), no.~16, 9635--9652. 

\bibitem{my-froz}
Kuznetsova M. Inverse problem for Sturm--Liouville operators with frozen argument on closed sets, Itogi Nauki Tekh., Ser. Sovrem. Mat. Prilozh., Temat. Obz. 208 (2022), 49--62.

\bibitem{bon-disc}
Bondarenko N.P. Finite-difference approximation of the inverse Sturm-Liouville problem with frozen argument, Appl. Math. Comput. 413
(2022), 126653.


\bibitem{AML} Kuznetsova M. Necessary and sufficient conditions for the spectra of the Sturm--Liouville operators with frozen argument, Applied Mathematics Letters 131 (2022), 108035.

\bibitem{DobHry} Dobosevych O., Hryniv R.  Reconstruction of differential operators with frozen argument, Axioms 11 (2022), no.1, article~24.

\bibitem{MN}  Buterin S.A. On the uniform stability of recovering sine-type functions with asymptotically separated zeros, Math. Notes 111 (2022), 343--355. 
\end{thebibliography}
\end{document}